\documentclass[leqno]{amsart}
       
\usepackage{amsmath,amsfonts,amssymb}
\usepackage{amsthm}
\usepackage{hyperref}

\def\R{{\mathbb R}}

\def\virgp{\raise 2pt\hbox{,}}

\def\({\left(}
\def\){\right)}
\def\<{\left\langle}
\def\>{\right\rangle}
\def\le{\leqslant}
\def\ge{\geqslant}

\def\Tend#1#2{\mathop{\longrightarrow}\limits_{#1\rightarrow#2}}

\def\d{{\partial}}

\def\g{\gamma}
\def\eps{\varepsilon}
\def\l{\lambda}

\def\si{{\sigma}}

\def\F{\mathcal F}
\def\M{\mathcal M}
\def\V{\mathcal V}
\def\conj{\mathcal C}
\def\O{\mathcal O}

\theoremstyle{plain}
\newtheorem{theorem}{Theorem}[section]
\newtheorem{lemma}[theorem]{Lemma}
\newtheorem{corollary}[theorem]{Corollary}
\newtheorem{proposition}[theorem]{Proposition}

\theoremstyle{definition}

\theoremstyle{remark}
\newtheorem{remark}[theorem]{Remark}
\newtheorem*{remark*}{Remark}

\numberwithin{equation}{section}

\begin{document}

\title[Wave operators for critical NLS]{On the wave operators for the
  critical nonlinear Schr\"odinger equation}   
\author[R. Carles]{R{\'e}mi Carles}
\address{D\'epartement de Math\'ematiques, UMR CNRS 5149\\CC~051\\
  Universit\'e Montpellier~2\\Place Eug\`ene Bataillon\\34095
  Montpellier cedex 5\\France\footnote{Present address: Wolfgang
  Pauli Institute, Universit\"at Wien, 
        Nordbergstr.~15, A-1090 Wien}}
\email{Remi.Carles@math.cnrs.fr}
\author[T. Ozawa]{Tohru Ozawa}
\address{Department of Mathematics\\ Hokkaido University\\ Sapporo
  060-810\\Japan}
\email{ozawa@math.sci.hokudai.ac.jp}
\begin{abstract}
 We prove that for the $L^2$-critical nonlinear Schr\"odinger
 equations, the wave operators and their inverse are related
 explicitly in terms of the Fourier transform. We discuss some
 consequences of this property. In the one-dimensional case, we show a
 precise similarity between the  $L^2$-critical nonlinear Schr\"odinger
 equation and a nonlinear Schr\"odinger
 equation of derivative type. 
\end{abstract}
\maketitle

\section{Introduction}
\label{sec:intro}
 We consider the defocusing, $L^2$-critical, nonlinear Schr\"odinger
 equation
\begin{equation}
  \label{eq:nlscrit}
  i\d_t u+\frac{1}{2}\Delta u = |u|^{4/n}u,\quad (t,x)\in \R\times
  \R^n. 
\end{equation}
We consider two types of initial data:
\begin{align}
  &\text{Asymptotic state: }U_0(-t)u(t)\big|_{t=\pm
  \infty}=u_\pm,\quad \text{where
  }U_0(t)=e^{i\frac{t}{2}\Delta}.\label{eq:CIasym} \\
&\text{Cauchy data at }t=0:\ u_{\mid t=0} = u_0.\label{eq:CIcauchy}
\end{align}
It is well known that for data $u_\pm,u_0\in \Sigma=H^1\cap\F(H^1)$, 
where
\begin{equation*}
  \F f(\xi)=\widehat
  f(\xi)=\frac{1}{(2\pi)^{n/2}}\int_{\R^n}f(x)e^{-ix\cdot \xi} dx,
\end{equation*}
\eqref{eq:nlscrit}--\eqref{eq:CIasym} has a unique, global, solution
$u\in C(\R;\Sigma)$ (\cite{GV79Scatt}, see also
\cite{CazCourant}). Its initial value $u_{\mid t=0}$ is the image of 
the asymptotic state under the action of the wave operator:
\begin{equation*}
  u_{\mid t=0}= W_\pm u_\pm. 
\end{equation*}
Similarly, \eqref{eq:nlscrit}--\eqref{eq:CIcauchy} possesses
asymptotic states:
\begin{equation*}
  \exists u_\pm\in \Sigma,\quad
  \|U_0(-t)u(t)-u_\pm\|_\Sigma \Tend t {\pm \infty} 0: \quad
  u_\pm=W_\pm^{-1}u_0.
\end{equation*}
Global well-posedness properties show that the wave operators are
homeomorphisms on $\Sigma$. Besides this point, very few properties of
these operators are known. The main result of this paper (proved in
\S\ref{sec:main}) shows that
the wave operators and their inverses are easily related in terms of
the Fourier transform:
\begin{theorem}\label{theo:main}
  Let $n\ge 1$. The following identity holds on $\Sigma$:
  \begin{equation}\label{eq:formula1}
    \F\circ W_\pm^{-1}=W_\mp\circ \F.
  \end{equation}
In particular, if $\conj$ denotes the conjugation $f\mapsto \overline
f$, then we have:
\begin{equation}\label{eq:formula2}
 W_\pm^{-1} = \(\conj  \F\)^{-1}  W_\pm \(\conj \F\). 
\end{equation}
\end{theorem}
Using continuity properties of the flow map associated to
\eqref{eq:nlscrit}, we infer the following result in \S\ref{sec:cor1}:
\begin{corollary}\label{cor:1}
 The result of Theorem~\ref{theo:main} still holds when $\Sigma$ is
 replaced
 \begin{itemize}
\item Either by $\F(H^1)$,
 \item Or by a neighborhood of the origin in $L^2(\R^n)$, for
 \eqref{eq:nlscrit} as well as for its focusing counterpart,
$  i\d_t u+\frac{1}{2}\Delta u = -|u|^{4/n}u$,
\item Or by $L^2_r(\R^n)$ for $n\ge 3$, the set of \emph{radial},
  square integrable functions.
 \end{itemize}
\end{corollary}
\begin{remark}
  The usual conjecture on \eqref{eq:nlscrit} implies that the result of
Theorem~\ref{theo:main} is expected to remain valid when $\Sigma$ is
replaced by $L^2(\R^n)$ (but not for the focusing counterpart of
\eqref{eq:nlscrit}, for which finite time blow-up may occur in $H^1$). 
\end{remark}
\begin{remark}
  So far, the existence of wave operators on
  $\F(H^1)$ is not known. Similarly, asymptotic completeness in $H^1$
  remains an open problem. Theorem~\ref{theo:main} shows that the fact
  that these two problems are simultaneously open is not merely a
  technical point: they are exactly related by
  \eqref{eq:formula1}. This aspect is also reminiscent of the main
  result in \cite{BlueColliander}. 
\end{remark}
Using the asymptotic expansion of the wave operators near the origin,
we prove in \S\ref{sec:cor2} (with an extension in
Appendix~\ref{sec:sub}): 
\begin{corollary}\label{cor:2}
Let $n\ge 1$. For every $\phi\in L^2(\R^n)$, we have:
\begin{equation*}
  \int_0^{\pm\infty} e^{it\frac{|x|^2}{2}}\F\(
  |U_0(t)\phi|^{4/n}U_0(t)\phi\) dt =
\int_0^{\pm \infty} U_0(t)\(
  |U_0(-t)\widehat\phi|^{4/n}U_0(-t)\widehat\phi\) dt.
\end{equation*}
\end{corollary}
\bigbreak

Finally, in space dimension $n=1$, we relate the wave operators for
\eqref{eq:nlscrit} with the wave operators for the nonlinear
Schr\"odinger equation of derivative type
\begin{equation}\label{eq:nlsder}
  i\d_t \psi +\frac{1}{2}\d_x^2 \psi = i \l\d_x\(|\psi|^2\)\psi,\quad
  \l\in \R. 
\end{equation}
This equation appears as a model to study the nonlinear
self-modulation for the Benjamin-Ono equation \cite{Tanaka82}. 
For a more general nonlinear
Schr\"odinger equation of derivative type (see
e.g. \cite{KatayamaTsutsumi94,Tsutsumi94} for the Cauchy problem
related to similar equations),
\begin{equation*}
  i\d_t \psi +\frac{1}{2}\d_x^2 \psi = i \l |\psi|^2\d_x \psi+ i\mu
  \psi^2\d_x \overline \psi,
\end{equation*}
it is proved in \cite{OzawaIndiana} that a short range scattering
theory is available for $\l,\mu\in\R$ if and only if $\l=\mu$: we
recover \eqref{eq:nlsder}. This is
apparently the only cubic, gauge 
invariant nonlinearity in space dimension one, for which a short range
scattering theory is available. More precisely, for
\eqref{eq:nlsder}--\eqref{eq:CIasym}, the wave operators
$\Omega_\pm: u_\pm\mapsto u(0)$ are well defined from $X_\eps$ to
$H^2(\R)$, where 
\begin{equation*}
  X_\eps = \{ \phi\in H^4\cap \F(H^4)\quad;\quad \left\| (1+\xi^2)\widehat
  \phi\right\|_{L^\infty}<\eps\},
\end{equation*}
and $\eps>0$ is sufficiently small. The following result shows that the
nonlinearity in \eqref{eq:nlsder} should be thought of as the quintic
case \eqref{eq:nlscrit}. This result goes in the same spirit as the
approach followed in \cite{OzawaTsutsumi98}. 
\begin{theorem}\label{theo:main2}
  Let $\l\in\R$. Consider the quintic, focusing or defocusing, equation
  \begin{equation}\label{eq:nlsmu}
   i\d_t u+\frac{1}{2}\d_x^2 u = \frac{\l^2}{2}|u|^4 u,\quad
  (t,x)\in \R\times 
  \R, \quad \mu\in \R,
  \end{equation}
with associated wave operators $W_\pm(\mu)$ for small $L^2$ data. For
$\phi\in L^2(\R)$, define 
\begin{equation*}
 ( N_\pm^\l \phi)(x) = \phi(x) \exp\( \pm i\l \int_{-\infty}^x
 |\phi(y)|^2dy \).  
\end{equation*}
$\bullet$ If $\psi$ solves \eqref{eq:nlsder}, then $N_-^\l(\psi)$ solves
\eqref{eq:nlsmu}.\\
$\bullet$ If $u$ solves \eqref{eq:nlsmu}, then
$N_+^\l(u)$ solves \eqref{eq:nlsder}.\\
$\bullet$ The following identity holds when all terms are well-defined:
\begin{align*}
 &\F \circ \Omega_\pm^{-1}= N_-^\l\circ \F\circ W_\pm^{-1} \circ N_+^\l=
 \(N_+^\l\)^{-1}\circ \F\circ W_\pm^{-1} \circ N_+^\l .\\
& \Omega_\pm\circ \F^{-1}= 
 \(N_+^\l\)^{-1}\circ W_\pm\circ \F^{-1} \circ N_+^\l .
\end{align*}
\end{theorem}
This result is checked by elementary computations, so we leave out its
proof. 
\section{Proof of Theorem~\ref{theo:main}}
\label{sec:main}

The proof of Theorem~\ref{theo:main} relies on a series of
lemmas, which are stated, and proved, in a slightly different fashion in
\cite{TsutsumiSigma}. Introduce the transform $\Psi$ acting on
function of $(t,x)$ as:
\begin{equation}
  \label{eq:pseudo}
  \(\Psi u\)(t,x) = \frac{1}{(it)^{n/2}}e^{i\frac{|x|^2}{2t}}
  u\(\frac{-1}{t}\virgp \frac{x}{t}\),\quad \text{for }t\not =0. 
\end{equation}
\begin{lemma}\label{lem:1}
  For $n\ge 1$ and $\phi\in L^2(\R^n)$, we have:
  \begin{equation*}
    \lim_{t\to\pm \infty}\left\| U_0(t)\F^{-1}\phi(\cdot) -
    (\Psi \phi)(t,\cdot)\right\|_{L^2}=0.
  \end{equation*}
\end{lemma}
\begin{proof}
  We recall the standard decomposition of the free group, for $t\not =0$:
  \begin{equation*}
    U_0(t) = \M_t D_t \F \M_t,
  \end{equation*}
where $\M_t$ is the multiplication by $e^{i|x|^2/(2t)}$, and $D_t$ is
the dilation operator
\begin{equation*}
  D_t\phi (x)=\frac{1}{(it)^{n/2}}\phi \(\frac{x}{t}\). 
\end{equation*}
Noting that $\Psi \phi = \M D\phi$, Plancherel formula yields:
\begin{equation*}
  \left\| U_0(t)\F^{-1}\phi(\cdot) -
    (\Psi \phi)(t,\cdot)\right\|_{L^2} = \left\|
    \(M_t-1\)\F^{-1}\phi(\cdot)\right\|_{L^2}.
\end{equation*}
Since $|M_t(x)-1|\lesssim |x|/\sqrt t$, the lemma follows for $\phi
\in H^1(\R^n)$. By density, we infer the result for  $\phi\in L^2(\R^n)$. 
\end{proof}
\begin{lemma}\label{lem:2}
  Let $v=\Psi u$. Suppose that there exist $\psi_\pm\in L^2(\R^n)$
  such that
  \begin{equation*}
    \|v(t)-\psi_\pm\|_{L^2}\Tend t {\pm 0} 0.
  \end{equation*}
Then $u$ has asymptotic states in $L^2$: 
\begin{equation*}
  \|U_0(-t)u(t) - \F^{-1}R \psi_\mp\|_{L^2}\Tend t {\pm \infty} 0,
\end{equation*}
where $R$ stands for the symmetry with respect to the origin,
$(R\phi)(x) = \phi(-x)$. 
\end{lemma}
\begin{proof}
  We note that $\Psi$ is almost an involution: $\Psi^2 =R$. Therefore,
  $u = \Psi Rv$:
  \begin{align*}
    U_0(-t)u(t) - \F^{-1}R \psi_\mp &= U_0(-t)\Psi R v\(\frac{-1}{t}\)
    - \F^{-1}R \psi_\mp \\
&= U_0(-t)\Psi R \(v\(\frac{-1}{t}\)-\psi_\mp\) + \(U_0(-t)\Psi-
    \F^{-1}\)R \psi_\mp . 
  \end{align*}
Taking the $L^2$ norm, we infer:
\begin{align*}
  \left\|U_0(-t)u(t) - \F^{-1}R \psi_\mp\right\|_{L^2}\le
  \left\|v\(\frac{-1}{t}\)-\psi_\mp\right\|_{L^2} + \|\Psi R \psi_\mp
  - U_0(t)\F^{-1}R \psi_\mp\|_{L^2}.
\end{align*}
The first term of the right-hand side goes to zero as $t\to \pm\infty$
by assumption. The second term goes to zero by Lemma~\ref{lem:1}.
\end{proof}
\begin{lemma}\label{lem:3}
  Let $v=\Psi u$. Suppose that $u\in C([-T,T];L^2)$ for some $T>0$,
  and $u_{\mid t=0}= u_0\in L^2(\R^n)$. Then
  \begin{equation*}
    \left\| U_0(-t)v(t)- \F^{-1}u_0\right\|_{L^2}\Tend t {\pm \infty}0.
  \end{equation*}
\end{lemma}
\begin{proof}
  Since $U_0(-t)=U_0(t)^{-1}$, we have
  \begin{equation*}
    U_0(-t)v(t) = \M_t^{-1}\F^{-1}D_t^{-1}\M_t^{-1} v(t) =
    \M_t^{-1}\F^{-1} u\(\frac{-1}{t}\).
  \end{equation*}
Therefore, 
\begin{equation*}
  \left\| U_0(-t)v(t) -  \F^{-1}u_0\right\|_{L^2} \le
  \left\|u\(\frac{-1}{t}\) -u_0\right\|_{L^2}+\left\|
\(\M_{-t}-1\) \F^{-1} u_0\right\|_{L^2} . 
\end{equation*}
The first term of the right-hand side goes to zero as $t\to \pm \infty$
by assumption. So does the second, by the standard argument
recalled in the proof of Lemma~\ref{lem:1}.
\end{proof}
\begin{proof}[Proof of Theorem~\ref{theo:main}] Let $u_0\in \Sigma$:
  there exists a unique solution $u\in C(\R;\Sigma)$ to
  \eqref{eq:nlscrit}--\eqref{eq:CIcauchy}. Set $v=\Psi u$. Because of
  the conformal invariance for \eqref{eq:nlscrit}, $v$ solves the same
  equation as $u$, for $t\not =0$:
  \begin{equation*}
    i\d_t v+\frac{1}{2}\Delta v = |v|^{4/n}v,\quad (t,x)\in
  \R\setminus\{0\}\times  \R^n. 
  \end{equation*}
 Lemma~\ref{lem:3} shows that 
 \begin{equation*}
   \left\| U_0(-t)v(t)- \F^{-1}u_0\right\|_{L^2}\Tend t {\pm \infty}0.
 \end{equation*}
Let $w_\pm$ denote the solutions to the scattering problems:
\begin{equation*}
  i\d_t w_\pm+\frac{1}{2}\Delta w_\pm = |w_\pm|^{4/n}w_\pm \quad
  ;\quad U_0(-t)w_\pm (t)\big|_{t=0}= \F^{-1}u_0.
\end{equation*}
By uniqueness for \eqref{eq:nlscrit}--\eqref{eq:CIasym}, we see that 
\begin{equation*}
  v(t,x)=\left\{
    \begin{aligned}
      w_-(t,x)& \quad \text{for }t<0,\\
      w_+(t,x)& \quad \text{for }t>0.
    \end{aligned}
\right.
\end{equation*}
In particular, 
\begin{equation*}
  \left\| v(t)-w_\pm(0)\right\|_{L^2} \Tend t {\pm 0}0. 
\end{equation*}
From Lemma~\ref{lem:2}, $u$ has asymptotics states, given by:
\begin{equation*}
  \|U_0(-t)u(t) - \F^{-1}Rw_\mp\|_{L^2}\Tend t {\pm \infty} 0,
\end{equation*}
that is, $u_\pm = \F^{-1}Rw_\mp$. We infer:
\begin{equation*}
  \F \circ W_\pm^{-1}u_0 = \F u_\pm = Rw_\mp= R W_\mp \F^{-1}u_0. 
\end{equation*}
Since \eqref{eq:nlscrit} is invariant by $R$, $R W_\mp \F^{-1}u_0 =
W_\mp R \F^{-1}u_0 = W_\mp \F u_0$. This yields \eqref{eq:formula1}.
The identity \eqref{eq:formula2} follows from \eqref{eq:formula1} and
from the identity
\begin{equation*}
  W_\pm = \conj \circ W_\mp\circ\conj ,
\end{equation*}
which was noticed in \cite{CW92} (see also \cite{CazCourant}). 
\end{proof}

\section{Proof of Corollary~\ref{cor:1}}
\label{sec:cor1}
The first case follows by density, since $W_\pm$ are defined and continuous on
$H^1(\R^n)$ \cite{GV85} (see also \cite{Ginibre} for a simplified
presentation), and since $W_\pm^{-1}$ are defined and continuous on
$\F(H^1)$ \cite{GOV94,GV79Scatt,TsutsumiSigma}.  
\smallbreak

For the second case, existence of wave operators, their asymptotic
 completeness, and continuity properties, were proved by T.~Cazenave and
F.~Weissler \cite{CW89}. We note that Corollary~\ref{cor:1} can be
proved in this case like Theorem~\ref{theo:main}, provided that we work in a
sufficiently small neighborhood of the origin in $L^2(\R^n)$.
\smallbreak

The last case follows from the recent paper
by T.~Tao, M.~Visan and X.~Zhang \cite{TVZradial}.
The proof of Corollary~\ref{cor:1} then relies on 
asymptotic completeness  (in the same space), along with continuous
dependence upon the 
initial data. For $n\ge 3$, let $X= L^2_r(\R^n)$; $X$ is invariant
under the action 
of the Fourier transform. For $\phi\in X$, let $\phi_j$ be a sequence
in $\Sigma$, converging to $\phi$ in $X$. Define $u^\pm_j$ as the
solutions to:
\begin{equation*}
  i\d_t u^\pm_j+\frac{1}{2}\Delta u^\pm_j =|u^\pm_j|^{4/n}u^\pm_j\quad
  ;\quad U_0(-t)u^\pm_j(t)\big|_{t=\pm \infty} = \widehat\phi_j. 
\end{equation*}
There exists $u^\pm_{0j}=u^\pm_j(0)= W_\pm  \widehat\phi_j \in
\Sigma$. Since $u^\pm_{0j} = \F W_\mp^{-1} \phi_j$ from
Theorem~\ref{theo:main}, the results in \cite{CW89,TVZradial} imply
that there exists $u^\pm_0\in X$ such 
that $\|u^\pm_{0j}-u^\pm_0\|_{L^2}\to 0$ as $j\to \infty$. Let $u^\pm$
solve 
\begin{equation*}
  i\d_t u^\pm+\frac{1}{2}\Delta u^\pm =|u^\pm|^{4/n}u^\pm\quad
  ;\quad u^\pm_{\mid t=0} = u^\pm_0. 
\end{equation*}
We have
\begin{align*}
  \|U_0(-t)u^\pm(t)- \widehat \phi\|_{L^2}\le &\ \left\|U_0(-t)\(u^\pm(t)
  - u^\pm_j(t)\)\right\|_{L^2} + \|U_0(-t)u^\pm_j(t)- \widehat
  \phi_j\|_{L^2}\\
& + \|\phi_j-\phi\|_{L^2}. 
\end{align*}
The global well-posedness for \eqref{eq:nlscrit} in $X$ implies 
\begin{align*}
  \limsup_{t\to \pm\infty }\|U_0(-t)u^\pm(t)- \widehat \phi\|_{L^2}\le
  & F\( \|u^\pm_0-u^\pm_{0j}\|_{L^2}\) + \|\phi_j-\phi\|_{L^2},
\end{align*}
where $F$ is a continuous function such that $F(0)=0$. Finally, by
letting $j\to \infty$, we see that $u^\pm$ solves
\begin{equation*}
  i\d_t u^\pm+\frac{1}{2}\Delta u^\pm =|u^\pm|^{4/n}u^\pm\quad
  ;\quad U_0(-t)u^\pm(t)\big|_{\mid t=\pm \infty} \widehat \phi. 
\end{equation*}
Let $\V$ be a neighborhood of $\phi$ in $L^2$. From \cite{CW89},
we see by Strichartz estimates and a 
bootstrap argument that the problem
\eqref{eq:nlscrit}--\eqref{eq:CIasym} is well-posed in
$L^\infty(]-\infty,-T];\V)$ (we consider only the minus sign for
simplicity) for 
some $T>0$ possibly depending on $\V$. By uniqueness, we infer 
\begin{equation*}
  \exists W_\pm \widehat \phi = u^\pm_0. 
\end{equation*}
Since under our assumptions, $W_\pm^{-1}$ are homeomorphisms on $X$,
we also have:
\begin{equation*}
 u^\pm_0 =\lim_{j\to \infty}  u^\pm_{0j} = \lim_{j\to \infty}\F
 W_\mp^{-1} \phi_j = \F
 W_\mp^{-1}\lim_{j\to \infty}\phi_j =\F
 W_\mp^{-1}\phi,
\end{equation*}
hence Corollary~\ref{cor:1}.

\section{Proof of Corollary~\ref{cor:2}}
\label{sec:cor2}
Corollary~\ref{cor:2} is a consequence of Theorem~\ref{theo:main} and
of the asymptotic expansion of the wave operators near the origin in
$L^2$: 
\begin{proposition}\label{prop:smallwave}
  Let $n\ge 1$ and $\phi\in L^2(\R^n)$. Then for $\eps>0$ sufficiently
  small $W_\pm (\eps^{n/4}\phi)$ and $W_\pm^{-1} (\eps^{n/4}\phi)$ are
  well defined in $L^2(\R^n)$, and, as $\eps\to 0$:
  \begin{align*}
    &W_\pm \(\eps^{\frac{n}{4}}\phi\) = \eps^{\frac{n}{4}}\phi \mp
 i\eps^{1+\frac{n}{4}}\int_0^{\pm 
 \infty} U_0(-t)\( |U_0(t)\phi|^{4/n}U_0(t)\phi\)dt +
 \O\(\eps^{2+\frac{4}{n}}\), \\   
 &W_\pm^{-1} \(\eps^{\frac{n}{4}}\phi\) = \eps^{\frac{n}{4}}\phi \pm
 i\eps^{1+\frac{n}{4}}\int_0^{\pm 
 \infty} U_0(-t)\( |U_0(t)\phi|^{4/n}U_0(t)\phi\)dt +
 \O\(\eps^{2+\frac{4}{n}}\).   
  \end{align*}
\end{proposition}
\begin{proof}
  The proof follows from the same perturbative analysis as in
  \cite{PG96} (see also \cite{CaWigner} for the nonlinear Schr\"odinger
  equation). First, it follows from \cite{CW89} that $W_\pm
  (\eps^{n/4}\phi)$ and $W_\pm^{-1} (\eps^{n/4}\phi)$ are
  well defined in $L^2(\R^n)$ for $\eps>0$ sufficiently small. 
\smallbreak

We prove the asymptotic formula for the minus sign, since the proof of
the formula for the plus sign is similar. Consider $u^\eps$ solving:
\begin{equation*}
  i\d_t u^\eps+\frac{1}{2}\Delta u^\eps = |u^\eps|^{4/n}u^\eps \quad
  ;\quad U_0(-t)u^\eps(t)\big|_{t=-\infty}=\eps^{n/4}\phi. 
\end{equation*}
Plugging an expansion of the form $u^\eps = \eps^{n/4}(\varphi_0 +
\eps\varphi_1 + \eps r^\eps)$ into the above equation, and ordering in
powers of $\eps$, it is natural to impose the following conditions:
\begin{itemize}
\item Leading order: $\O(\eps^{n/4})$.
  \begin{equation*}
    i\d_t \varphi_0 + \frac{1}{2}\Delta \varphi_0 = 0\quad ;\quad
    U_0(-t)\varphi_0(t)\big|_{t=-\infty}=\phi .
  \end{equation*}
\item First corrector: $\O(\eps^{1+n/4})$.
  \begin{equation*}
   i\d_t \varphi_1 + \frac{1}{2}\Delta \varphi_1 =
    |\varphi_0|^{4/n}\varphi_0\quad ;\quad 
    U_0(-t)\varphi_1(t)\big|_{t=-\infty}=0 . 
  \end{equation*}
\end{itemize}
The first equation yields
\begin{equation*}
  \varphi_0(t)= U_0(t)\phi.
\end{equation*}
From the second equation, we have:
\begin{equation*}
  \varphi_1(t) = -i\int_{-\infty}^t U_0(t-s)\(
  |\varphi_0(s)|^{4/n}\varphi_0(s)\) ds. 
\end{equation*}
We also have:
\begin{equation*}
  i\d_t r^\eps +\frac{1}{2}\Delta r^\eps = G\( \varphi_0 +
\eps\varphi_1 + \eps r^\eps \) - G(\varphi_0)\quad ;\quad 
    U_0(-t)r^\eps(t)\big|_{t=-\infty}=0,
\end{equation*}
where $G(z)=|z|^{4/n}z$. Let $\g =2+4/n$, and denote $L^r_{t,x}=
L^r(]-\infty,-t]\times\R^n)$. Strichartz and H\"older estimates yield
\begin{align*}
  \|r^\eps\|_{L^\g_{t,x}}&\lesssim \left\| \(|\varphi_0|^{4/n} +
  |\eps\varphi_1|^{4/n}+ |\eps r^\eps|^{4/n} \) \eps \(|\varphi_1|
  +|r^\eps|\)\right\|_{L^{\g'}_{t,x}} \\
& \lesssim \(
  \|\varphi_0\|_{L^\g_{t,x}}^{4/n}+
\|\eps\varphi_1\|_{L^\g_{t,x}}^{4/n}+
\|\eps r^\eps\|_{L^\g_{t,x}}^{4/n}\)\( \eps \|\varphi_1\|_{L^\g_{t,x}}
  +\|\eps  r^\eps\|_{L^\g_{t,x}}  \)\\
& \lesssim  \(\|\phi\|_{L^2}^{4/n}+\|\eps
  r^\eps\|_{L^\g_{t,x}}^{4/n}\)\( \eps\|\phi\|_{L^2}
  +\| \eps r^\eps\|_{L^\g_{t,x}}  \).\\
 & \lesssim \eps \|\phi\|_{L^2}^{1+4/n}+
\| \eps r^\eps\|_{L^\g_{t,x}}^{1+4/n}.
\end{align*}
A bootstrap argument shows that for $0<\eps\ll 1$, 
$r^\eps \in L^\g (\R\times\R^n)$, and
\begin{equation*}
  \|r^\eps\|_{L^\g(\R\times\R^n)}\lesssim \eps.
\end{equation*}
Using Strichartz estimates again, we infer:
\begin{equation*}
  \|r^\eps\|_{L^\infty(\R;L^2(\R^n))}\lesssim \eps.
\end{equation*}
Considering $u^\eps$ at time $t=0$ yields the first part of the
proposition. The second part can be proven in the same way, but can
also be inferred from the first part \emph{via} Neumann series, since
$W_\pm$ are small perturbations of the identity near the origin. 
\end{proof}
Now Corollary~\ref{cor:2} follows from Corollary~\ref{cor:1} and
Proposition~\ref{prop:smallwave}, where we identify the terms of order
$\eps^{1+4/n}$. 
\begin{remark}
  Considering the asymptotic expansion of the wave operators and their
  inverse to higher order would yield other formulae, similar to
  Corollary~\ref{cor:2}. We have not written them, for they are more
  intricate (they involve several integrations in time), and we do not
  know if they can be of some interest. 
\end{remark}

\appendix

\section{Sub-critical case}
\label{sec:sub}
In this appendix, we consider more generally the nonlinear
Schr\"odinger equation
\begin{equation}
  \label{eq:nls}
  i\d_t u+\frac{1}{2}\Delta u = |u|^{2\si}u,\quad (t,x)\in \R\times
  \R^n,
\end{equation}
in the sub-critical case $\si<2/n$. Following the approach to prove
Corollary~\ref{cor:2}, we have:
\begin{proposition}
  Let $\si<2/n$, with 
  \begin{itemize}
  \item $\si >1/n$ if $n\le 2$.
\item $\si >2/(n+2)$ if $n\ge 2$. 
  \end{itemize}
Then the following identities hold for every $\phi \in \Sigma$:
 \begin{align*}
&  \int_0^{\pm\infty} e^{it\frac{|x|^2}{2}}\F\(
  |U_0(t)\phi|^{2\si}U_0(t)\phi\) dt =
\int_0^{\pm \infty} |t|^{n\si -2}U_0(t)\(
  |U_0(-t)\widehat\phi|^{2\si}U_0(-t)\widehat\phi\) dt,\\
&\int_0^{\pm\infty} |t|^{n\si -2} e^{it\frac{|x|^2}{2}}\F\(
  |U_0(t)\phi|^{2\si}U_0(t)\phi\) dt =\int_0^{\pm \infty} U_0(t)\(
  |U_0(-t)\widehat\phi|^{2\si}U_0(-t)\widehat\phi\) dt.
\end{align*} 
\end{proposition}
\begin{proof}[Sketch of the proof]
  Let $u$ solving \eqref{eq:nls}. Then $v=\Psi u$ solves
\begin{equation}
  \label{eq:nls2}
  i\d_t v+\frac{1}{2}\Delta v = |t|^{n\si -2}|v|^{2\si}v,\quad
  (t,x)\in \R\setminus\{0\}\times 
  \R^n.
\end{equation}
It follows from \cite{CW92} (see also \cite{CazCourant}) that wave
operators exist, are continuous and invertible, near the origin in
$\Sigma$, both for \eqref{eq:nls} and \eqref{eq:nls2}. We can then
mimic the proof of Theorem~\ref{theo:main}, with the remark that in
Theorem~\ref{theo:main}, the operators $W_\pm^{-1}$ on the left-hand side are
associated to $u$, while the operators  $W_\mp$ on the right-hand side
are associated to $v$.  
\smallbreak

Adapting Proposition~\ref{prop:smallwave} to the cases of
\eqref{eq:nls} and \eqref{eq:nls2} proceeds along the same lines as
the estimates in \cite{CW92}. This yields the first identity in the
above proposition. 
\smallbreak

For the second, we simply notice that $\Psi^2 = R$, so that we can
exchange the roles of $u$ and $v$. 
\end{proof}

\providecommand{\bysame}{\leavevmode\hbox to3em{\hrulefill}\thinspace}
\providecommand{\href}[2]{#2}


\begin{thebibliography}{GOV94}

\bibitem[BC06]{BlueColliander}
P.~Blue and J.~Colliander, \emph{Global well-posedness in {S}obolev space
  implies global existence for weighted ${L}^2$ initial data for
  ${L}^2$-critical {NLS}}, Commun. Pure Appl. Anal. \textbf{5} (2006), no.~4,
  691--708.

\bibitem[Car01]{CaWigner}
R.~Carles, \emph{Remarques sur les mesures de {W}igner}, C. {R}. {A}cad. {S}ci.
  {P}aris, t. 332, {S}\'erie {I} \textbf{332} (2001), no.~11, 981--984.

\bibitem[Caz03]{CazCourant}
T.~Cazenave, \emph{Semilinear {S}chr\"odinger equations}, Courant Lecture Notes
  in Mathematics, vol.~10, New York University Courant Institute of
  Mathematical Sciences, New York, 2003.

\bibitem[CW89]{CW89}
T.~Cazenave and F.~Weissler, \emph{Some remarks on the nonlinear
  {S}chr\"odinger equation in the critical case}, Lect. Notes in Math., vol.
  1394, Springer-Verlag, Berlin, 1989, pp.~18--29.

\bibitem[CW92]{CW92}
\bysame, \emph{Rapidly decaying solutions of the nonlinear {S}chr\"odinger
  equation}, Comm. Math. Phys. \textbf{147} (1992), 75--100.

\bibitem[G{\'e}r96]{PG96}
P.~G{\'e}rard, \emph{Oscillations and concentration effects in semilinear
  dispersive wave equations}, J. Funct. Anal. \textbf{141} (1996), no.~1,
  60--98.

\bibitem[Gin97]{Ginibre}
J.~Ginibre, \emph{An introduction to nonlinear {S}chr\"odinger equations},
  Nonlinear waves (Sapporo, 1995) (R.~Agemi, Y.~Giga, and T.~Ozawa, eds.),
  GAKUTO International Series, Math. Sciences and Appl., Gakk\={o}tosho, Tokyo,
  1997, pp.~85--133.

\bibitem[GOV94]{GOV94}
J.~Ginibre, T.~Ozawa, and G.~Velo, \emph{On the existence of the wave operators
  for a class of nonlinear {S}chr\"odinger equations}, {A}nn. {IHP} ({P}hysique
  {T}h\'eorique) \textbf{60} (1994), 211--239.

\bibitem[GV79]{GV79Scatt}
J.~Ginibre and G.~Velo, \emph{On a class of nonlinear {S}chr\"odinger
  equations. {II} {S}cattering theory, general case}, J. Funct. Anal.
  \textbf{32} (1979), 33--71.

\bibitem[GV85]{GV85}
\bysame, \emph{Scattering theory in the energy space for a class of nonlinear
  {S}chr\"odinger equations}, J. Math. Pures Appl. (9) \textbf{64} (1985),
  no.~4, 363--401.

\bibitem[KT94]{KatayamaTsutsumi94}
S.~Katayama and Y.~Tsutsumi, \emph{Global existence of solutions for nonlinear
  {S}chr\"odinger equations in one space dimension}, Comm. Partial Differential
  Equations \textbf{19} (1994), no.~11-12, 1971--1997.

\bibitem[OT98]{OzawaTsutsumi98}
T.~Ozawa and Y.~Tsutsumi, \emph{Space-time estimates for null gauge forms and
  nonlinear {S}chr\"odinger equations}, Differential Integral Equations
  \textbf{11} (1998), no.~2, 201--222.

\bibitem[Oza96]{OzawaIndiana}
T.~Ozawa, \emph{On the nonlinear {S}chr\"odinger equations of derivative type},
  Indiana Univ. Math. J. \textbf{45} (1996), no.~1, 137--163.

\bibitem[Tan82]{Tanaka82}
M.~Tanaka, \emph{Nonlinear self-modulation problem of the {B}enjamin-{O}no
  equation}, J. Phys. Soc. Japan \textbf{51} (1982), no.~8, 2686--2692.

\bibitem[Tsu85]{TsutsumiSigma}
Y.~Tsutsumi, \emph{Scattering problem for nonlinear {S}chr\"odinger equations},
  Ann. Inst. H. Poincar\'e Phys. Th\'eor. \textbf{43} (1985), no.~3, 321--347.

\bibitem[Tsu94]{Tsutsumi94}
\bysame, \emph{The null gauge condition and the one-dimensional nonlinear
  {S}chr\"odinger equation with cubic nonlinearity}, Indiana Univ. Math. J.
  \textbf{43} (1994), no.~1, 241--254.

\bibitem[TVZ]{TVZradial}
T.~Tao, M.~Visan, and X.~Zhang, \emph{Global well-posedness and scattering for
  the mass-critical nonlinear {S}chr\"odinger equation for radial data in high
  dimensions}, preprint, {\tt math.AP/0609692}.

\end{thebibliography}
\end{document}